\documentclass[final]{svjour3PP}

\usepackage{amsfonts,amsmath,amssymb}

\usepackage{mathptmx}
\usepackage{cite}

\AtBeginDocument{%
  \setlength{\oddsidemargin}{\dimexpr(\paperwidth-\textwidth)/2-1in}%
  \setlength{\evensidemargin}{\oddsidemargin}%
  \setlength{\topmargin}{%
    \dimexpr(\paperheight-\textheight)/2-\headheight-\headsep-1in}%
}

\def\0{\ensuremath{\mathbf{0}}}

\usepackage{graphicx}
\usepackage{a4,amssymb,amsmath}
\def\be{\begin{equation}}
\def\bea{\begin{eqnarray*}}
\def\ee{\end{equation}}
\def\eea{\end{eqnarray*}}
\def\ba{\begin{array}}
\def\ea{\end{array}}
\def\bi{\begin{itemize}}
\def\ei{\end{itemize}}
\newtheorem{theo}{Theorem}
\newtheorem{lem}{Lemma}
\newtheorem{pro}{Proposition}

\def\RR{\mathbb{R}}
\def\NN{\mathbb{N}}
\def\TT{\mathbb{T}}

\def\teps{\tilde{\epsilon}}

\def\gL{\Lambda}

\def\mA{\mathcal{A}}
\def\mD{\mathcal{D}}

\title{Schwarz Iterative Methods: Infinite Space Splittings}
\author{Michael Griebel \and Peter Oswald}
\institute{M. Griebel\at Institute for Numerical Simulation, Universit\"at Bonn, Wegelerstr. 6, 53115 Bonn, Germany, and Fraunhofer Institute for Algorithms and Scientific Computing (SCAI), Schloss Birlinghoven, 53754 Sankt Augustin, Germany\\
Corresponding author, tel.: +49-228-733437, fax: +49-228-737527,\\ 
\email{griebel@ins.uni-bonn.de}
\and P. Oswald \at Jacobs University Bremen, Research I, 600 Campus Ring 1, D-28759 Bremen\\
\email{p.oswald@jacobs-university.de, agp.oswald@gmail.com}}

\titlerunning{Schwarz Iterative Methods: Infinite Space Splittings}
\authorrunning{M. Griebel \and P. Oswald}

\date{}
\begin{document}
\maketitle
\begin{abstract}
We prove the convergence of greedy and randomized versions of Schwarz iterative methods for solving linear elliptic variational problems
based on infinite space splittings of a Hilbert space. For the greedy case, we show a squared error decay rate of $O((m+1)^{-1})$ for elements of an approximation space 
$\mA_1$ related to the underlying splitting.
For the randomized case, we show an expected squared error decay rate of $O((m+1)^{-1})$ on a class
$\mA_{\infty}^{\pi}\subset \mA_1$ depending on the probability distribution. %Moreover, we discuss further results for more special situations.

\end{abstract}

\section{Introduction}\label{sec1}
The aim of this paper is to extend convergence results for greedy and randomized versions of multiplicative Schwarz methods
for solving elliptic variational problems in Hilbert spaces from the case of finite space splittings \cite{GrOs2011,OsZh2014} to the case of infinite space splittings. 
Let $V$ be a separable real or complex Hilbert space with scalar product $(\cdot,\cdot)$, let 
$a(\cdot,\cdot)$ be a continuous and coercive Hermitian form on $V$, and let $F$ be a bounded linear functional on $V$.
Note that $a(\cdot,\cdot)$ induces a spectrally equivalent scalar product on $V$, in the sequel we write
$V_a$ to indicate that we consider $V$ with this new scalar product, and use the notation $\|\cdot\|_a$ for the 
induced norm.
Then the variational problem
\bi
\item[(A)] Find $u\in V$ such that
$$
a(u,v)=F(v)\qquad\forall v\in V,
$$
\ei
possesses a unique solution, and is equivalent to the quadratic minimization problem
\bi
\item[(B)] Find the minimizer $u\in V$ of the quadratic functional
$$
\Phi(u):=\frac12a(u,u)-F(u).
$$
\ei
We treat the problem in the form (A), by turning it into an infinite linear system using space splittings
as described next. The equivalent formulation (B) provides the link to convex optimization, where algorithms similar
to the ones considered here are known and investigated under the name block-coordinate descent methods.

For the separable Hilbert space $V_a$, we consider {\it space splittings} generated by families
of Hilbert spaces $V_{a_i}$ (with scalar product $a_i(\cdot,\cdot)$ and norm $\|\cdot\|_{a_i}$) and bounded linear operators
$R_i:\,V_{a_i}\to V_a$, $i\in I$, such that the span of the subspaces $R_iV_{a_i}\subset V_a$ is dense in $V_a$. Here, the index set $I$ can be finite (${I}=\{1,2,\ldots,N\}$), or countable (${I}=\NN$). These conditions on a space splitting are silently assumed throughout
this paper. We call a space splitting {\it stable}, and write
\be\label{SSS}
V_a=\sum_{i\in I} R_iV_{a_i} := \left\{v=\sum_{i\in I} R_iv_i\,:\;v_i\in V_{a_i},\;i\in I\right\},
\ee
if
\be\label{SSS1}
0< \lambda_{\min}:= \inf_{u\in V_a} \frac{a(u,u)}{\||u\||^2}\;\le\;\lambda_{\max}:= \sup_{u\in V_a} \frac{a(u,u)}{\||u\||^2}<\infty,
\ee
where
$$
\||u\||^2 := \inf_{v_i\in V_{a_i}: u=\sum_{i\in {I}} R_iv_i} \sum_{i\in {I}} a_i(v_i,v_i).
$$
Not every infinite space splitting is stable, in particular, in (\ref{SSS}) it is assumed that every element in $V_a$ possesses at least one converging series expansion
with respect to $\{R_iV_{a_i}\}$ which follows from (\ref{SSS1}) and the assumed density of span($\{R_iV_{a_i}\}$) in $V_a$. The constants $\lambda_{\min}$ and $\lambda_{\max}$ are called lower and upper stability constants, and
$\kappa:=B/A$ is called the condition of the stable space splitting (\ref{SSS}), respectively. A prominent
case of stable space splittings are frames and fusion frames, see \cite{Ch,CaKu,PO,PO3}.
In all these definitions, we allow for redundancy, i.e., $R_iV_{a_i}\cap R_jV_{a_j} = \{0\}$ is not required for $i\neq j$, and we do not assume that the $V_{a_i}$
are closed subspaces of $V_a$.

For the setup of Schwarz iterative methods we need
 to define operators $T_i:\,V_a\to V_{a_i}$ via the variational problems
\be\label{subproblem}
a_i(T_iv,v_i)=a(v,R_iv_i)  \qquad \forall\;v_i\in V_{a_i},
\ee
to be solved for given $v\in V_a$ on the spaces $V_{a_i}$, $i\in {I}$. Evaluating $T_iv$ is equivalent to solving
a variational problem in $V_{a_i}$, and it is silently assumed that this is easier than solving the original
problem (A). This stems from the fact that $V_{a_i}$ has typically much smaller dimension than $V_a$ and/or the Hermitian form $a_i(\cdot,\cdot)$
leads to a linear system with better spectral properties or simpler structure.
If the underlying space splitting is finite then, using these $T_i$, analogs of the classical Jacobi-Richardson and Gauss-Seidel-SOR iterations,
called additive and multiplicative Schwarz methods with respect to (stable) space splittings can be defined and investigated,
pretty much along the lines of the standard methods, see \cite{Xu,Gr,PO,GrOs,PO4}.

Here we formulate a generic version of the multiplicative (also called sequential or asynchronous) Schwarz method with relaxation suitable
for the case of infinite space splittings.
Choose an initial approximation $u^{(0)}$, and repeat the following steps for $m=0,1,\ldots$ until a stopping criterion is met:
\begin{itemize}
\item[1.] {\bf Subproblem pick and solution}: Given the current $u^{(m)}$ and an index set $I_m\subset I$, choose an index $i=i_m\in I_m$
(according to some rule to be specified),
and compute the partial residual $r^{(m)}_i:=T_ie^{(m)}$, where
$e^{(m)}:= u-u^{(m)}$. Although $u$ is unknown, this can be done since the right-hand side in the
corresponding subproblem (\ref{subproblem}) reads
$$
a(e^{(m)},R_iv_i)=F(R_iv_i)-a(u^{(m)},R_iv_i),
$$
and does not depend on knowledge about $u$.
\item[2.] {\bf Linear update}: Determine relaxation parameters $\alpha_m\ge 0$ and $\omega_m$ (according to some rule to be specified), and set
\be\label{updateLin}
u^{(m+1)}=\alpha_m u^{(m)}+\omega_m R_{i}r^{(m)}_{i}.
\ee
\end{itemize}
So far, 
this is a theoretical algorithm since executing Steps 1 and 2 is not feasible without further specification
and assumptions. In the case of finite splittings, the Schwarz iterative method figures also under the name alternating directions method (ADM), see \cite{Ga2004} for references.

As to Step 1, we need to specify the rule for picking the next index $i=i_m$. There are at least three standard versions to be considered:
\bi
\item {\it Deterministic orderings}. In this case, we choose an index sequence $\{i_m\}$ beforehand.
In the case of finite splittings, the default orderings are cyclic ($i_{rN+k-1}=k$) or
symmetric-cyclic ($i_{2rN+k-1}=k,\;i_{2rN+N+k-1}=N+1-k$) for $k=1,\ldots,N$, $r=0,1$) which corresponds to the classical SOR and SSOR methods, respectively. A naive deterministic ordering for infinite space splittings would be
to choose $\{1;1,2;1,2,3;1,2,3,4;\ldots\}$. For finite splittings, convergence is known for cyclic orderings from the ADM theory (compare, e.g., \cite{Ga2004}),
see also the convergence rate estimates for Schwarz iterative methods in \cite{GrOs} and, more recently, for coordinate descent methods and convex optimization 
problems \cite{BeTe2013}.
\item {\it Greedy orderings}. The idea goes back to Gauss and Seidel, and was popularized by Southwell \cite{So,So1}
(the corresponding algorithms for finite splittings are often called Gauss-Southwell methods).
Here the decision for the next index  $i_m $ depends on the current iterate $u^{(m)}$, and aims at
maximizing the error reduction in the next step. For instance, we can require $i_m\in I_m$ to satisfy
\be\label{greedySSS}
a_{i_m}(r^{(m)}_{i_m},r^{(m)}_{i_m})\ge \beta_m^2 \sup_{i\in I_m} a_i(r^{(m)}_i,r^{(m)}_i),
\ee
where $\beta_m\in (0,1]$ is called weakness parameter. 
This approach is expensive, as it involves the computation of multiple partial residuals, at least approximately, just to pick the next index.
Most of the research on quantitative convergence results for greedy orderings and infinite splittings (see \cite{Te2011} for an overview) is devoted to the case $I_m=I=\NN$, where finding an $i_m$ that satisfies (\ref{greedySSS}) in a numerically feasible way can be guaranteed only under additional assumptions. In practice, 
one would prefer working with dynamically growing but finite index sets $I_m\subset I$.
In the case of finite splittings, algorithms with greedy orderings have been analyzed in the more general setting of convex optimization methods, see e.g. \cite{Ts,LuTs,Ts2001} for early results in this direction, for a short proof in the case of problem (A), see \cite{GrOs2011}.
\item {\it Random orderings}. Choose a sequence of discrete probability distributions 
$$\pi^{(m)}=\{\pi_i^{(m)}\ge 0\}_{i\in I},$$ 
and pick $i=i_m\in I$ randomly according to $\pi^{(m)}$, $m=0,1,\ldots$. 
Even in the case of finite splittings, a theoretical analysis of such algorithms has been started only recently but it revealed that they are competitive with the best (often unknown)  deterministic orderings, and numerically much cheaper than greedy orderings. We refer to \cite{LeLe,GrOs2011,OsZh2014,StVe} for the setting of the present paper (quadratic minimization as in (B)), and to  \cite{Ma2013,Ne2012,RiTa2012,TaRiGo2013} for recent convergence results on block coordinate
descent methods for large-scale convex optimization problems.
\ei
Certainly, there are many more variants to explore. In this paper, we concentrate on greedy and
random orderings for infinite splittings ($I=\NN$).

As to Step 2, many options have been discussed in the literature, especially in connection with greedy orderings, see
\cite{Te,Te2011} for an overview and references. 
\bi
\item The simplest algorithms result if we fix both parameters $\alpha_m=1$, $\omega_m=\omega$ independently of $m$,
and assume some normalization condition for the $R_i$. Then we arrive at analogs of the algorithms discussed in the theory of greedy methods
under the names weak (WGA) and weak relaxed (WRGA) greedy algorithms. Their counterparts for $\beta=1$ are called pure (PGA) and relaxed (RGA) greedy algorithms, respectively.
\item More generally, one can try to find $\alpha_m\ge 0, \omega_m$ simultaneously  by minimizing the new error term
\be\label{errornext}
\|u-u^{(m+1)}\|_a=\min_{\alpha\ge 0,\omega} \|u-\alpha u^{(m)}-\omega R_{i}r_{i}^{(m)}\|_a,
\ee
with respect to $\alpha$ and $\omega$. This is equivalent to solving a two-dimensional quadratic minimization problem.
Various restrictions have been proposed, e.g., $\alpha+\omega=1$ is a popular choice, see the early work \cite{Jo,Ba} on relaxed greedy algorithms, and \cite{Te2012a,DeTe2014} in a slightly more general setting. 
\item In this paper, we consider a variant with fixed parameter sequence
\be\label{alphachoice}
\alpha_m:=1-(m+2)^{-1} \in (0,1), \qquad m\ge 0,
\ee
and with $\omega_m$ determined by minimizing the error term (\ref{errornext}) with respect to $\omega$. More explicitly,
\be\label{omegachoice}
\omega_m = \frac{a(u-\alpha_m u^{(m)},R_{i}r_{i}^{(m)})}{\|R_{i}r_{i}^{(m)}\|_{a}^2}
=\frac{\alpha_m a_{i}(r_{i}^{(m)},r_{i}^{(m)})+ \bar{\alpha}_m F(R_{i}r^{(m)}_{i})}{\|R_{i}r_{i}^{(m)}\|_{a}^2}, 
\ee
where $\bar{\alpha}_m:=(1-\alpha_m)=(m+2)^{-1}$, $m\ge 0$. In the theory of greedy algorithms, this version is labeled GAWR, see \cite{BCDD,Te2}.
\ei
Further modifications are listed and investigated in \cite{Te2,Te2011}, for extensions to convex optimization problems on Hilbert and Banach spaces 
see the recent papers \cite{Zh2003,Te2012a,DeTe2014,NgPe2014}. We also mention that instead of finding an optimal approximation $u^{(m+1)}$ only from 
span($\{R_ir_i^{(m)},u^{(m)}\}$) as done in (\ref{errornext}), one could include earlier approximations $u^{(k)}$, $k<m$, into the local search in Step 2.
Orthogonal matching pursuit is a relatively expensive extension of this type. A less expensive version motivated by the conjugate gradient 
method would be to find $u^{(m+1)}$ from span($\{R_ir_i^{(m)},u^{(m)}-u^{(m-1)}\}$).

The main contributions of this paper can be summarized as follows.
For the greedy case, we restrict ourselves in (\ref{greedySSS}) to 
$$
\beta_m=\beta \in (0,1],\qquad I_m=\NN, \qquad  m\ge 0,
$$
and give a convergence proof for the above specified 
theoretical algorithm. This proof is, in the case of the Hilbert space setting, a modification of the approach
used in \cite{BCDD,Te2}. In particular, we show the convergence rate 
$$
\epsilon_m:=\|u-u^{(m)}\|_a =\mathrm{O}((m+1)^{-1/2}),\qquad m\to \infty,
$$
for $u$ from the class $\mA_1$  which will be defined below. This class appears naturally in all investigations on greedy algorithms, and the
exponent $1/2$ in the convergence rate estimate is known to be optimal in our considered situation of general space splittings \cite{DeTe}. 

For the case of random picks
with a fixed probability distribution $\pi^{(m)}=\pi>0$, we prove a similar estimate for the expected error decay,
$$
\tilde{\epsilon}_m := E(\epsilon_m^2)^{1/2} =\mathrm{O}((m+1)^{-1/2}),\qquad m\to \infty,
$$
for $u$ from a smaller class $\mA_{\infty}^{\pi}\subset \mA_1$. 
To the best of our knowledge, 
this is the first general convergence result for randomized Schwarz iterative methods in the case of infinite splittings. Using an approximation and density argument, we show convergence in expectation (without guaranteed rate) also 
for arbitrary $u\in V_a$ and for sequences $\pi^{(m)}$ that converge to a fixed probability distribution $\pi>0$ in $\ell^1$.

Although mathematically not difficult, we emphasize that our approach via space splittings (rather than expansions with respect to dictionaries $\mD\subset V$ and updates along one-dimensional search directions) covers block-iterative methods, auxiliary space techniques, and outer approximation schemes which may lead to a broader applicability of our theoretical findings.

The remainder of this paper is organized as follows:
In Section \ref{sec2} we present our convergence results. To this end, we first introduce approximation spaces $\mA_q^\gamma$ associated with a given infinite  space splitting, and give some preparatory lemmata. Then, we prove the main theorems on convergence estimates for greedy and  randomized Schwarz iterations.
In Section \ref{sec3} we discuss some further results and consequences.
%Finally, we give some concluding remarks in Section \ref{sec4}.

\section{Convergence Results}\label{sec2}

\subsection{Approximation spaces related to space splittings}\label{sec21}
Throughout this section, set $I=\NN$, and fix the families $\{V_{a_i}\}_{i\in\NN}$ of auxiliary Hilbert spaces and $\{R_i:\,V_{a_i}\to V\}_{i\in\NN}$ of bounded linear operators. Furthermore, let them be such that the span of the linear subspaces $R_iV_{a_i}$ is dense in $V_a$.
We assume uniform boundedness of the operators $R_i$, i.e., there exists a constant $\gL$ such that
\be\label{RBound}
\|R_iv_i\|_a^2 =a(R_iv_i,R_iv_i) \le \gL^2 a_i(v_i,v_i)=\gL^2 \|v_i\|^2_{a_i},\qquad v_i\in V_{a_i},\quad i\in \NN.
\ee
In theory, this can always be achieved by rescaling either $R_i$ or the auxiliary Hermitian forms $a_i$.  For the practical application this is however irrelevant as the minimization with respect to $\omega$ in the update Step 2
of the algorithms automatically takes care of this. Unless stated otherwise, we will {\em not} assume that (\ref{SSS}) is a stable space splitting for $V_a$.

For any non-negative weight sequence $\gamma=\{\gamma_i\}_{i\in\NN}\ge 0$, and $0<q\le \infty$, introduce approximation spaces $\mA^\gamma_q\equiv \mA^\gamma_q(\{V_{a_i},R_i;V_a\})$ as follows:
For $u\in \mathrm{span}(\{R_iV_{a_i}\}_{i\in\mathrm{supp}(\gamma)})$, define 
$$
|u|_{\mA^\gamma_q} = \inf \left\{\|\{\gamma_i^{-1}\|v_i\|_{a_i}\}_{i\in I}\|_{\ell^q}:\, v_i\in V_{a_i}\; u = \sum_{i\in I} R_iv_i 
\quad(\mbox{with finite\;}I \subset \mathrm{supp}(\gamma) )\right\},
$$
and introduce $\mA^\gamma_q$ as the completion with respect to this (quasi-)semi-norm. For the weight sequence $\gamma_i=1$, we drop the superscript $\gamma$ from the notation.
The cases we are most interested in are 
$$
\mA^{\pi}_\infty \subset \mA_1 \subset \left\{\ba{l} V_a,\\ \mA_q,\qquad 1<q\le 2, \ea\right. 
$$
where $\pi\ge 0$ is any given discrete probability distribution, i.e., $\pi_i\ge 0$ and $\sum_{i\in \NN} \pi_i=1$. The embeddings are continuous. If (\ref{SSS}) is a stable splitting, then obviously
$\mA_2=V_a$, and all spaces $\mA^{\pi}_\infty$ and $\mA_q$, $q<2$, are subspaces of $V_a$. They are also dense in $V_a$ (for $\mA^{\pi}_\infty$ under the assumption that $\mathrm{supp}(\pi)=\NN$).
For the case of one-dimensional subspaces of $V_a$ generated by a countable dictionary, these definitions are standard and have been instrumental for setting up a quantitative convergence theory for greedy algorithms with infinite dictionaries, see \cite{Te,Te1}. 

The following technical lemma is crucial for our convergence proofs below.
\begin{lem}\label{lem1}
For the underlying space splitting, assume (\ref{RBound}).
For any $e\in V_a$, denote $r_i=T_ie\in V_{a_i}$, $w_i= \|r_{i}\|_{a_{i}}^{-1} R_ir_i\in V_a$, $i\in \NN$.\\
a) If $i^\ast$ is such that $\|r_{i^\ast}\|_{a_{i^\ast}}\ge \beta \sup_{i\in \NN} \|r_{i}\|_{a_i}$ for some $0<\beta\le 1$, then, for any nontrivial $h\in \mA_1$, we have
\be\label{Agreedy}
\|r_{i^\ast}\|_{a_{i^\ast}} = a(e,w_{i^\ast}) \ge \frac{\beta}{\|h\|_{\mA_1}} a(e,h).
\ee
b) If $\pi$ is any discrete probability distribution, then, for any nontrivial $h\in \mA_\infty^{\pi}$, we have
\be\label{Arandom}
\sum_{i\in \NN} \pi_i\|r_{i}\|_{a_{i}} = \sum_{i\in \NN} \pi_i a(e,w_i) \ge \frac{1}{\|h\|_{\mA_\infty^{\pi}}} a(e,h).
\ee
\end{lem}

{\bf Proof}. Since $r_i=T_ie$ is the unique minimizer of the associated quadratic minimization problem, i.e.,
$$
\frac12 a_i(T_ie,T_ie) - a(e,R_iT_ie) \le \frac12 a_i(v_i,v_i) - a(e,R_iv_i)\qquad \forall\; v_i\in V_{a_i},
$$
for any index $i$, one concludes for any $v_i\in V_{a_i}$ with $\|v_i\|_{a_i} \le \|r_i\|_{a_i}$ that
$$
\|r_i\|_{a_i}^2=a(e,R_iT_ie)=\|r_i\|_{a_i} a(e,w_i) \ge a(e,R_iv_i).
$$
If $r_i\neq 0$, after dividing by $\|r_i\|_{a_i}$, this yields
$$
\|r_i\|_{a_i}=a(e,w_i)\ge a(e,R_iv_i) \qquad \forall v_i\in V_{a_i}\,:\; \|v_i\|_{a_i}\le 1.
$$
This inequality holds for $r_i=0$ as well, since in this case $0=a_i(T_ie,v_i)=a(e,R_iv_i)$ for all $v_i\in V_{a_i}$.
Since by
$$
\|\sum_{i\in \NN} R_iv_i\|_a \le \sum_{i\in \NN} \|R_iv_i\|_a\le \gL \sum_{i\in \NN} \|v_i\|_{a_i}
$$
convergence in $\mA_1$ obviously implies convergence 
in $V_a$, we can assume that for any $\epsilon>0$, there exists a finitely representable 
\be\label{hhref}
\tilde{h}=\sum_{i\in I'} c_iR_iv_i
\ee
(with finite $I'\subset \NN$ and $\|v_i\|_{a_i}\le 1$) such that 
$$
\gL^{-1}\|h-\tilde{h}\|_a\le \|h-\tilde{h}\|_{\mA^1} < \epsilon,\qquad \sum_{i\in I'} |c_i| \le \|h\|_{\mA_1} +\epsilon.
$$
For Part a), we thus have
$$
a(e,w_{i^\ast}) \ge  \beta \sup_i a(e,w_i) \ge \beta \sup_{\|v_i\|_{a_i}\le 1}  a(e,R_iv_i) \ge \frac{\beta }{\|h\|_{\mA_1}+\epsilon}a(e,\tilde{h}),
$$
and letting $\epsilon \to 0$ implies (\ref{Agreedy}).

Similarly, in Part b) we can choose $\tilde{h}$ of the form (\ref{hhref}) such that
$$
\gL^{-1}\|h-\tilde{h}\|_a\le \|h-\tilde{h}\|_{\mA_\infty^{\pi}} < \epsilon,\qquad |c_i| \le \pi_i(\|h\|_{\mA_\infty^{\pi}}+\epsilon).
$$
Then, by the same reasoning
$$
\sum_{i\in \NN} \pi_i a(e,w_i) \ge \sup_{\|v_i\|_{a_i}\le 1} \sum_{i\in \NN} \pi_i a(e,R_iv_i) = \sup_{\|v_i\|_{a_i}\le \pi_i} a(e,\sum_{i\in \NN} R_i v_i)
\ge \frac{1}{\|h\|_{\mA_\infty^{\pi}}+\epsilon} a(e,\tilde{h}).
$$
With $\epsilon\to 0$, we get (\ref{Arandom}) which finishes the proof of Lemma 1. \hfill$\Box$

Below, we will apply this lemma with $e=e^{(m)}:=u-u^{(m)}$,
i.e. the error after $m$ steps of the algorithm, and with $\beta=\beta_m$, i.e. the weakness parameter in the case of greedy orderings, while $\pi$ coincides
with the probability distribution used to create random orderings.

As another preparation, we formulate an auxiliary result for approximation in $\mA^{\pi}_\infty$ spaces that will allow us to work with variable probability distributions (see also Remark 4 in Section \ref{sec3}).
\begin{lem}\label{lem2}
Assume that $\pi>0$ is a fixed probability distribution with support $\NN$, and assume that $\pi^{(m)}\ge 0$ is a sequence of probability distributions that converges to $\pi$ in the $\ell^1$ norm. For the underlying space splitting, assume (\ref{RBound}). Then, for any given $h\in \mA^{\pi}_\infty$, there exists a sequence of finitely representable 
$$
h^{(m)} = \sum_{i\in I^{(m)}} R_iv_i, \qquad I^{(m)} \mbox{\; finite},
$$
such that for $m\ge 0$
\be\label{Hm}
\|h-h^{(m)}\|_{a} \le (1+3\gL)\|h\|_{\mA^{\pi}_\infty}\|\pi-\pi^{(m)}\|_{\ell^1},\qquad \|h^{(m)}\|_{{\mA^{\pi}_\infty}} \le 3\|h\|_{{\mA^{\pi}_\infty}}.
\ee
\end{lem}

{\bf Proof}. Since $h\in \mA^{\pi}_\infty$, for given $m\ge 0$ and $\delta>0$ there is a finite $I'$ and a 
$$
\bar{h}=\sum_{i\in I'} R_iv_i, \qquad \|v_i\|_{a_i} \le (1+\delta)\pi_i \|h\|_{\mA^{\pi}_\infty}, \quad i\in I',
$$
such that 
$$
\|h-\bar{h}\|_{a} \le \|h\|_{\mA^{\pi}_\infty}\|\pi-\pi^{(m)}\|_{\ell^1}.
$$
In the definition of $h^{(m)}$, set $I^{(m)}=\{i\in \NN: \pi^{(m)}_i> c\pi_i\}\cap I'$ with a constant $c$ to be fixed below. Then, obviously, we have
\bea
\|h-h^{(m)}\|_{a} &\le & \|h-\bar{h}\|_a +\sum_{i\in I'\backslash I^{(m)}} \|R_iv_i\|_a \le \|h-\bar{h}\|_a + \gL\sum_{i\in I^{(m)}} \|v_i\|_{a_i}\\
&\le & \|h\|_{\mA^{\pi}_\infty} \left(\|\pi-\pi^{(m)}\|_{\ell^1}+ (1+\delta)\gL \sum_{i\in I'\backslash I^{(m)}} \pi_i\right).
\eea
But for $i\in I'\backslash I^{(m)}$ we have $(1-c)\pi_i \le \pi_i-\pi^{(m)}_i$, thus
$$
\|h-h^{(m)}\|_{a} \le \|h\|_{\mA^{\pi}_\infty} \|\pi-\pi^{(m)}\|_{\ell^1} \left(1+\frac{(1+\delta)\gL}{1-c}\right).
$$
On the other hand, by the definition of $\bar{h}$ and $h^{(m)}$
$$
\|h^{(m)}\|_{\mA^{\pi^{(m)}}_\infty}\le \max_{i\in I^{(m)}} \frac{\|v_i\|_{a_i}}{\pi^{(m)}_i} \le
\|h\|_{\mA^{\pi}_\infty} \frac{1+\delta}{c}.
$$
Choosing $\delta=c=1/2$ gives then the statement of Lemma \ref{lem2}.\hfill $\Box$

\subsection{Convergence Estimates}\label{sec22}
As in the case of finite space splittings \cite{GrOs2011}, our convergence proof of the Schwarz iterative method for infinite space splittings
is based on the same error representation for both greedy and random orderings.
We therefore state the core estimates together in one place.

\begin{theo}\label{theo1}
Consider an infinite space splitting consisting of auxiliary Hilbert spaces $V_{a_i}$ and bounded linear operators $R_i:\,V_{a_i}\to V_a$, $i\in\NN$, 
such that span($\{R_iV_{a_i}\}_{i\in\NN}$) is dense in $V_a$ and (\ref{RBound}) holds. Furthermore, consider a Schwarz iterative method for the variational problem (A) with starting approximation $u^{(0)}=0$ and update rule
(\ref{updateLin}), where the parameters $\alpha_m$ and $\omega_m$ are specified by (\ref{alphachoice}) and (\ref{omegachoice}), respectively.\\
a) Assume that the update indices $i=i_m$ are chosen according to the greedy rule (\ref{greedySSS}) with $I_m=\NN$ and weakness parameter $\beta_m=\beta\in (0,1]$, $m\ge 0$.
If $u\in \mA_1$ then the squared error decay is given by
$$
\|u-u^{(m)}\|_a^2 \le 2(\|u\|_a^2 +(\gL/\beta)^2 \|u\|_{\mA_1}^2) (m+1)^{-1},\qquad m\ge 0.
$$
b) Assume that the update indices $i=i_m\in \NN$ are chosen randomly and independently according to a fixed discrete probability distribution $\pi$, $m\ge 0$.
If $u\in {\mA^{\pi}_\infty}$ then the expected squared error decay is given by
$$
E(\|u-u^{(m)}\|_a^2) \le 2(\|u\|_a^2 +\gL^2 \|u\|_{\mA^{\pi}_\infty}^2) (m+1)^{-1},\qquad m\ge 0.
$$
\end{theo}

{\bf Proof}. We derive a recursion for the (expected) squared error.
Suppose that $u^{(m)}$ is determined, and that the $i$-th subproblem solution $r_i^{(m)}:=T_ie^{(m)}$ is used for the update to $u^{(m+1)}$
according to (\ref{updateLin}).
Thus, we can write 
$$
e^{(m+1)}:=u-(\alpha_m u^{m} + \theta_{i,m} R_ir_i^{(m)}) = \alpha_m e^{(m)} +\bar{\alpha}_m(u-\xi_{i,m}w_i^{(m)}),
$$
where $\bar{\alpha}_m=1-\alpha_m$, and, in agreement with the previous subsection,  $w_i^{(m)}=\|r_i\|_{a_i}^{-1} R_ir_i^{(m)}$.  
The parameter $\xi_{i,m}=\theta_{i,m}/\bar{\alpha}_m$ is found by solving the minimization problem 
$$
\|e^{(m+1)}\|_a^2 = \min_\xi \|\alpha_m e^{(m)} +\bar{\alpha}_m(u-\xi w_i^{(m)})\|_a^2.
$$
Thus, for any $\xi$ and any chosen $i$, we have
$$
\|e^{(m+1)}\|_a^2 \le \alpha_m^2 \|e^{(m)}\|_a^2 +2\alpha_m\bar{\alpha}_m(a(e^{(m)},u)-\xi a(e^{(m)},w_i^{(m)})) + \bar{\alpha}_m^2\|u-\xi w_i^{(m)}\|_a^2.
$$
Using the triangle inequality and (\ref{RBound}), the norm in the last term can be bounded, independently of $i$, by
$$
\|u-\xi w_i^{(m)}\|_a^2\le 2(\|u\|_a^2+\xi^2 \| r_i^{(m)}\|_{a_i}^{-2}\| R_ir_i^{(m)}\|_{a}^2)\le 2(\|u\|_a^2+\xi^2 \gL^2).
$$

For dealing with the term $a(e^{(m)},u)-\xi a(e^{(m)},w_i^{(m)})$, we invoke Lemma \ref{lem1}. 
In the case of greedy orderings, we take $h=u\in \mA_1$ in Part a) of Lemma \ref{lem1} and choose $\xi=\beta^{-1}\|u\|_{\mA_1}$.
We then arrive at
\be\label{Egreedy}
\|e^{(m+1)}\|_a^2 \le \alpha_m^2 \|e^{(m)}\|_a^2 + 2\bar{\alpha}_m^2(\|u\|_a^2 + (\gL/\beta)^{2}\|u\|_{\mA^1}^2).
\ee
Denote $c_m:=(m+1)\|e^{(m)}\|_a^2$ and $M:=2(\|u\|_a^2 + (\gL/\beta)^{2}\|u\|_{\mA_1}^2)$. Substituting the concrete values of $\alpha_m=1-(m+2)^{-1}$
and $\bar{\alpha}_m=(m+2)^{-1}$ into (\ref{Egreedy}), we obtain
$$
c_{m+1} \le \alpha_m c_m + \bar{\alpha}_m M, \qquad m\ge 0.
$$
Since $c_0=\|e^{(0)}\|^2_a=\|u\|^2_a < M$, this implies $c_m <M$ for all $m\ge 0$, and proves the result in Part a) of Theorem \ref{theo1}.

In the case of random orderings, we now use Lemma \ref{lem1} b) for the given discrete probability distribution $\pi$ with $h=u\in {\mA^{\pi}_\infty}$.
This yields the following estimate for the conditional expectation of $\|e^{(m+1)}\|_a^2$ with respect to given $u^{(m)}$
valid for any $\xi\ge 0$:
\bea
 E(\|e^{(m+1)}\|_a^2\,|\, u^{(m)}) &\le& \alpha_m^2 \|e^{(m)}\|_a^2 
+2\alpha_m\bar{\alpha}_m (a(e^{(m)},u) -\xi \sum_i \pi_i a(e^{(m)},w_i^{(m)}))
\\ && \qquad\qquad\qquad\qquad+ 2\bar{\alpha}_m^2(\|u\|_a^2+\xi^2 \gL^2)\\
&\le& \alpha_m^2 \|e^{(m)}\|_a^2 +2\alpha_m\bar{\alpha}_m a(e^{(m)},u)(1-\xi \|u\|_{\mA^{\pi}_\infty}^{-1})
+ 2\bar{\alpha}_m^2(\|u\|_a^2+\xi^2 \gL^2).
\eea
Thus, fixing $\xi=\|u\|_{\mA^{\pi}_\infty}$ and taking the expectations with respect to $u^{(m)}$, we get
\be\label{Erandom}
E(\|e^{(m+1)}\|_a^2) \le \alpha_m^2 E(\|e^{(m)}\|_a^2) + 2\bar{\alpha}_m^2(\|u\|_a^2 + \gL^{2}\|u\|_{\mA^{\pi}_\infty}^2).
\ee
This gives Part b) of Theorem \ref{theo1} if we argue as before.\hfill $\Box$

\medskip
Using a density argument as in \cite{BCDD,Te2}, one can extend the estimate of Theorem \ref{theo1}, and show convergence for all $u\in V_a$. 

\begin{theo}\label{theo2}
Under the same assumptions as in Theorem \ref{theo1}, we have convergence and expected convergence $u^{(m)}\to u$ in $V_a$ for the greedy and random Schwarz iterative methods, respectively,
with no additional assumptions on the solution $u\in V_a$ of the variational problem (A). \\
More precisely, for the greedy version specified in 
Part a) of Theorem \ref{theo1} and any $h\in \mA_1$, we have
\be\label{Hgreedy}
\|u-u^{(m)}\|_a \le 2 \|u-h\|_a + \frac{\sqrt{8(\|u\|_a^2 + (\gL/\beta)^2\|h\|_{\mA_1}^2)}}{(m+1)^{1/2}},\qquad m\ge 0.
\ee
For the random version specified in 
Part b) of Theorem \ref{theo1} and any $h\in {\mA^{\pi}_\infty}$, we have
\be\label{Hrandom}
E(\|u-u^{(m)}\|_a^2)^{1/2} \le 2 \|u-h\|_a + \frac{\sqrt{8(\|u\|_a^2+\gL^2 \|h\|_{\mA^{\pi}_\infty}^2)}}{(m+1)^{1/2}},\qquad m\ge 0.
\ee
\end{theo}

{\bf Proof}. We start with the greedy case, take any
$h\in \mA^1$. Repeat the proof of Part a) of Theorem \ref{theo1}. When invoking Lemma \ref{lem1} a), use it with $h$ instead of $u$, and 
set $\xi=\beta^{-1}\|h\|_{\mA_1}$. Then 
$$
a(e^{(m)},u)-\xi a(e^{(m)},w_i^{(m)}) \le a(e^{(m)},u) - \xi \frac{\beta}{\|h\|_{\mA_1}} a(e^{(m)},h) = a(e^{(m)},u-h)\le \|e^{(m)}\|_a\|u-h\|_a,
$$
and (\ref{Egreedy}) can be replaced by
$$
\|e^{(m+1)}\|_a^2 \le \alpha_m^2 \|e^{(m)}\|_a^2 + 2\alpha_m\bar{\alpha}_m \|e^{(m)}\|_a\|u-h\|_a  +   2\bar{\alpha}_m^2(\|u\|_a^2 + (\gL/\beta)^{2}\|h\|_{\mA_1}^2).
$$
If $\|e^{(m+1)}\|_a > \alpha_m \|e^{(m)}\|_a$, then 
$$
\|e^{(m+1)}\|_a \le \frac{\|e^{(m+1)}\|_a^2}{\alpha_m \|e^{(m)}\|_a}.
$$
Substituting the previous estimate of $\|e^{(m+1)}\|_a^2$, this yields
\be\label{Egreedy1}
\|e^{(m+1)}\|_a \le \alpha_m \|e^{(m)}\|_a+ 2\bar{\alpha}_m \|u-h\|_a  +  \frac{\bar{\alpha}_m^2 M}{{\alpha}_m\|e^{(m)}\|_a},
\ee
where $M=2(\|u\|_a^2 + (\gL/\beta)^{2}\|h\|_{\mA_1}^2)$. If, alternatively,  $\|e^{(m+1)}\|_a \le \alpha_m \|e^{(m)}\|_a$,
then (\ref{Egreedy1}) holds trivially. The inequality (\ref{Egreedy1}) is complemented by
\be\label{Egreedy2}
\|e^{(m+1)}\|_a = \inf_\xi \|\alpha_m e^{(m)} +\bar{\alpha}_m(u-\xi w_i^{(m)})\|_a 
\le \alpha_m\|e^{(m)}\|_a +\bar{\alpha}_m\|u\|_a.
\ee

In the random case, we proceed similarly. For any $h\in {\mA^{\pi}_\infty}$, we apply Lemma 1 b) with $\xi=\|h\|_{\mA^{\pi}_\infty}\ge 0$. This shows
$$
a(e^{(m)},u)-\xi \sum_{i\in \NN} \pi_i a(e^{(m)},w_i^{(m)})\le a(e^{(m)},u-h)\le \|e^{(m)}\|_a\|u-h\|_a,
$$
and, instead of (\ref{Erandom}), we obtain
$$
E(\|e^{(m+1)}\|_a^2) \le \alpha_m^2 E(\|e^{(m)}\|_a^2) + 2\alpha_m\bar{\alpha}_m E(\|e^{(m)}\|_a)\|u-h\|_a  +   \bar{\alpha}_m^2\tilde{M}.
$$
where $\tilde{M}:=2(\|u\|_a^2 + \gL^{2}\|h\|_{\mA^{\pi}_\infty}^2)$. Using the notation $\teps_m:=E(\|e^{(m)}\|_a^2)^{1/2}$ and the obvious
inequality $E(\|e^{(m)}\|_a)\le \teps_m$, by the same reasoning as above,  we arrive at the following replacement for (\ref{Erandom}):
\be\label{Erandom1}
\teps_{m+1} \le \alpha_m \teps_m+ 2\bar{\alpha}_m \|u-h\|_a  +  \frac{\bar{\alpha}_m^2 \tilde{M}}{{\alpha}_m\teps_m}.
\ee
To obtain a complementary estimate analogous to (\ref{Egreedy2}), by definition of $u^{(m+1)}$ we can write
$$
\|e^{(m+1)}\|_a^2\le \|u-\alpha_m u^{(m)}\|_a^2\le \alpha_m^2\|e^{(m)}\|_a^2 + 2\alpha_m\bar{\alpha}_m \| u\|_a\|e^{(m)}\|_a + \bar{\alpha}_m^2 \|u\|_a^2.
$$
Then we take expectations on both sides, use again $E(\|e^{(m)}\|_a)\le \teps_m$, and obtain
\be\label{Erandom2}
\teps_{m+1} \le (\alpha_m^2 \teps_m^2 + 2\alpha_m\bar{\alpha}_m \| u\|_aE(\|e^{(m)}\|_a)+\bar{\alpha}_m^2 \| u\|_a^2)^{1/2}\le \alpha_m \teps_m + \bar{\alpha}_m \| u\|_a.
\ee

Up to different constants $M$ and $\tilde{M}$, the recursive inequalites (\ref{Egreedy1}-\ref{Egreedy2}) for the sequence
$\{\epsilon_m\}$ and (\ref{Erandom1}-\ref{Erandom2}) for $\{\teps_m\}$ are identical.
Therefore, it is enough to consider the random case. Set $b_m=(m+1)^{1/2}\tilde{M}^{-1/2}(\teps_m-2\|u-h\|_a)$,
$m\ge 0$. Then, a quick calculation shows that (\ref{Erandom2}) turns into
\be\label{E2}
b_{m+1} \le \alpha_m^{1/2} b_m + \frac{B}{(m+2)^{1/2}},\qquad m\ge 0,
\ee
where $B:=\|u\|_a \tilde{M}^{-1/2}$, while under the assumption $b_m>0$ the inequality (\ref{Erandom1}) implies 
\be\label{E1}
b_{m+1} \le \alpha_m^{1/2} (b_m+ \frac{1}{(m+1)b_m}).
\ee
A similar recursive system of inequalities has been considered in \cite{BCDD,Te2}. 
%The inequalities (\ref{E2}) and (\ref{E3}) imply
%\be\label{Best}
%b_m\le \max(b_0,\frac{B}{\sqrt{2}}+ \sqrt{1+\sqrt{2}}),\qquad m\ge 0.
%\ee
%Indeed, if $b_m\ge \sqrt{1+\sqrt{2}}$ then (\ref{E1}) shows 
%$$
%b_{m+1}-b_m\le \frac{\alpha_m^{1/2}}{b_m(m+1)}(1-(m+1) (\alpha_m^{-1/2}-1)b_m^2)\le \frac{\alpha_m^{1/2}}{b_m(m+1)}(1-\frac{2(1+\sqrt{2})\sqrt{m+1}}{\sqrt{m+1}+\sqrt{m+2}}\le 0,
%$$
%i.e., $b_{m+1}\le b_m$. On the other hand, if $b_m\le \sqrt{1+\sqrt{2}}$ then (\ref{E2}) implies
%$$ 
%b_{m+1}\le \alpha_m^{1/2} \sqrt{1+\sqrt{2}} + \frac{B}{(m+2)^{1/2}} \le \sqrt{1+\sqrt{2}} + \frac{B}{\sqrt{2}}.
%$$
% These two estimates imply (\ref{Best})
Lemma \ref{lem3} stated below implies
$b_m\le 2$ for all $m\ge 0$ (note that $b_0\le B\le 1/\sqrt{2}$). 
%Since in our case $b_0\le B \le 1/\sqrt{2}$, we conclude from (\ref{Best}) that
This yields
$$
\teps_m\le 2\|u-h\|_a + \frac{2\sqrt{\tilde{M}}}{(m+1)^{1/2}}, \qquad m\ge 0,
$$
and proves (\ref{Hrandom}). The estimate (\ref{Hgreedy}) is derived in complete analogy.

To show convergence for arbitrary $u\in V_a$,
for given $\epsilon >0$, choose $h\in \mA_1$ by the density of $\mA_1$ in $V_a$ such that $\|u-h\|_a<\epsilon/3$. Then, with this $h$ fixed, the second term in the right-hand side of (\ref{Hgreedy}) will become $<\epsilon/3$ for all large enough $m$ as well. This proves convergence for the greedy version. An analogous argument shows
$E(\|u-u^{(m)}\|_a^2)\to 0$ as $m\to \infty$ for the random case. Theorem \ref{theo2} is fully established.\hfill $\Box$

\medskip
For the convenience of the reader, we conclude this section with the short proof of the boundedness of sequences $\{b_m\}_{m\ge 0}$ 
satisfying the recursion (\ref{E2}-\ref{E1}) used in the proof of Theorem \ref{theo2}.

\begin{lem} \label{lem3} Suppose, a sequence $\{b_m\}_{m\ge 0}$ satisfies the inequalities (\ref{E2}) and (\ref{E1}), where $\alpha_m=(m+1)/(m+2)$ and $B>0$.
Fix a constant $A\ge B/\sqrt{2} +\sqrt{2}$. Then $b_0\le A$ implies $b_m\le A$ for all $m\ge 1$. In particular, if $B\le 1/\sqrt{2}$ one can choose $A=2$.
\end{lem}

{\bf Proof}. We use induction in $m$. Assume $b_m\le A$. For a value $t=t_m>0$ to be fixed below, we consider two cases. If $b_m\le t$, by (\ref{E2}) we have
$$
b_{m+1} \le \alpha_m^{1/2} t + \frac{B}{(m+2)^{1/2}} \le A,
$$
if
\be\label{C2}
t\le \alpha_m^{-1/2}( A - \frac{B}{(m+2)^{1/2}}).
\ee
On the other hand, if $t\le b_m\le A$ we use (\ref{E1}) which gives
$$
b_{m+1} \le \alpha_m^{1/2}( A + \frac{1}{(m+1)t} )\le A,
$$
if
\be\label{C1}
t\ge \frac{\alpha_m^{1/2}}{(1-\alpha_m^{1/2})(m+1) A} = \frac{(m+1)^{1/2}+(m+2)^{1/2}}{(m+1)^{1/2}A}=\alpha_m^{-1/2}\frac{(m+1)^{1/2}+(m+2)^{1/2}}{(m+2)^{1/2}A}.
\ee
It is easy to see that the choice $t=t_m:=2\alpha_m^{-1/2}/A>0$ satisfies both (\ref{C2}) and (\ref{C1}) if $2/A\le A - B/\sqrt{2}$. The latter follows from the
assumption on $A$ which implies $A\ge A-B/\sqrt{2} \ge \sqrt{2}\ge 2/A$. This finishes the induction step, and proves Lemma \ref{lem3}.  \hfill $\Box$

\section{Further Results and Discussion}\label{sec3}
{\bf Remark 1.}  Our results for the expected error decay for random orderings imply immediately estimates in probability. 
Using the Markov-Chebyshev inequality, under the assumptions of Theorem \ref{theo1} b), we get  
$$
\mathrm{P}\left(\|u-u^{(m)}\|_a^2 \le  \epsilon^2\right) \ge 1-\frac{8(\|u\|_a^2 +\gL^2 \|u\|_{\mA^{\pi}_\infty}^2)}{ (m+1)\epsilon^2},\qquad m\ge 0,
$$
for any error threshold $\epsilon>0$, or, equivalently,
$$
\mathrm{P}\left(\|u-u^{(m)}\|_a^2 \le \frac{8(\|u\|_a^2 +\gL^2 \|u\|_{\mA^{\pi}_\infty}^2)}{(m+1)\delta}\right) \ge 1-\delta,\qquad m\ge 0,
$$
for any confidence level $\delta$. An investigation of the variance or other higher-order moments of the squared error that could lead to improved
estimates has not been undertaken yet. Numerical experiments with randomized Schwarz iterations for finite splittings \cite{GrOs2011,OsZh2014} suggest that the 
variance is reasonably small in practice.

{\bf Remark 2.}
If in addition to our assumptions on space splittings we assume that (\ref{SSS}) is a stable space splitting of $V_a$, i.e., if $\mA_2=V_a$
holds, then, for $u\in \mA_q$, $1<q<2$ and the greedy version of the Schwarz iterative method specified in Part a) of Theorem \ref{theo1}, we have the error decay rate
\be\label{Aqerror}
\|u-u^{(m)}\|_a\le \bar{C} (m+1)^{1/2-1/q} \|u\|_{\mA_q}, \qquad m\ge 0,
\ee
where $\bar{C}$ is some absolute constant depending on $\beta$ and
the upper stability constant $\lambda_{\max}$ only. This can be established by an interpolation argument along the lines 
%
%
%To establish (\ref{Aqerror}), take $u\in \mA^q$ and consider a representation
%$$
%u=\sum_i R_iv_i = (u-h)+h,\qquad \sum_i \|v_i\|_{a_i}^q < 2\|u\|_{\mA_q}^q,
%$$
%where
%$$
%h:=\sum_{i\in I_1} R_iv_i, \qquad   I_1:=\{i:\,\|v_i\|_{a_i} \ge (m+1)^{-1/q}\|u\|_{\mA^q}\}
%$$
%(such a representation exists by the definition of $\mA^q$ and the assumption $V_a=\mA^2$). Observe that with this choice of $h$, we have
%\bea
%\|h\|_{\mA^1} &\le& \sum_{i\in I_1} \|v_i\|_{a_i}
%\le ((m+1)^{-1/q}\|u\|_{\mA^q})^{-(q-1)}\sum_{i\in I_1} \|v_i\|_{a_i}^q\\
%&\le& 2((m+1)^{-1/q}\|u\|_{\mA^q})^{-(q-1)}\|u\|_{\mA^q}^q
%=2(m+1)^{1-1/q}\|u\|_{\mA^q}<\infty,
%\eea
%and
%\bea
%\|u-h\|_{\mA^2}^2 &\le& \sum_{i\not\in I_1} \|v_i\|_{a_i}^2\le ((m+1)^{-1/q}\|u\|_{\mA^q})^{(2-q)}\sum_{i\not\in I_1} \|v_i\|_{a_i}^q\\
%&\le& 2((m+1)^{-1/q}\|u\|_{\mA^q})^{(2-q)}\|u\|_{\mA^q}^q
%=2(m+1)^{2(1/2-1/q)}\|u\|_{\mA^q}^2<\infty.
%\eea
%Substitution into (\ref{Hgreedy}) gives the statement,
%if one takes into account that
%\bea
%\|u-h\|_a&\le& \sqrt{\lambda_{\max}}\|u-h\|_{\mA^2},\\
%\|u\|_a &\le& \sqrt{\lambda_{\max}}\|u\|_{\mA^2}\le C(\|u-h\|_{\mA^2}+\|h\|_{\mA^1}).
%\eea
%Note that by the same token $\gL\le \lambda_{\max}$ is guaranteed.
of \cite{BCDD,Te2}, where the authors consider the special case of splittings into one-dimensional subspaces
$V_{a_i}=\{\lambda_i \psi_i:\;\lambda_i\in \RR\}$ of $V_a$ induced by a dictionary
$\mathcal{D}=\{\psi_i\}$ of unit norm elements $\psi_i\in V_a$ such that its span is dense in $V_a$. For this case, the estimate (\ref{Aqerror}) may be replaced by
a similar statement for $u\in \mathcal{B}_q$, where $\mathcal{B}_q$, $1<q<2$, is obtained by real interpolation
for the pair $(\mathcal{A}_1,V_a)$. In the special case, when
$\mathcal{D}$ is a frame in $V_a$ and thus $V_a=\mathcal{A}_2$, the scales $\mathcal{A}_q$ and $\mathcal{B}_q$
coincide for $1\le q<2$, in general, they are different. 

%It looks surprising that the convergence bound (\ref{Aqerror}) only relies on $\lambda_{\max}$. However,
%$\lambda_{\min}$ enters (with a negative power) if one wants to bound the $\mathcal{A}^q$ norm of $u$
%appearing in (\ref{Aqerror}) by the energy norm $\|u\|_a$.

\smallskip
{\bf Remark 3}. In \cite{Te}, weaker error estimates for other greedy algorithms such as PGA and WGA
can be found. We believe that their proofs carry over to the setting based on space splittings
without difficulties. E.g., for the PGA with $\alpha_m=\beta=1$, we expect
$$\|u-u^{(m)}\|_a \le Cm^{-1/6}\|u\|_{\mathcal{A}_1},\qquad m\ge 1,\qquad u\in \mathcal{A}_1,
$$
to hold, see \cite{DeTe} for the PGA in the dictionary case. Whether the exponent $1/6$ can be increased to $1/2$ under the assumption that
(\ref{SSS}) is a stable space splitting is an open problem, even when the space splitting
comes from a frame. Slightly better exponents are possible for PGA and WGA, see \cite{Te,Te2011}.

\smallskip
{\bf Remark 4}. Theorems \ref{theo1} and \ref{theo2} provide convergence guarantees under theoretical assumptions that look still questionable
from a practical point of view: The question of rounding errors is not addressed, for the greedy version the condition (\ref{greedySSS}) needs to
be checked for an infinite index set $I_m=\NN$, while in the random Schwarz iterative method drawing the next index $i=i_m$ according to a (rather general) discrete probability distribution $\pi$ defined on $\NN$ seems inconvenient as well. 
For greedy algorithms based on dictionaries, there are partial results in this direction \cite{Te2005,BCDD,DeTe2014} which can be adapted to the case of space splittings considered here.

We concentrate on the randomized version. When combined with Lemma \ref{lem2}, the estimation techniques leading to Theorem \ref{theo2} give the following result
which, in particular, allows us to work with finitely supported probability distributions $\pi^{(m)}$ that converge to a desired
$\pi>0$ sufficiently fast, without sacrificing convergence speed.
\begin{pro}\label{pro3} 
Assume that the indices $i=i_m$ in the random Schwarz iteration 
are chosen using discrete probability distributions $\pi^{(m)}\ge 0$ such that
\be\label{PiError}
\|\pi^{(m)}-\pi\|_{\ell^1} \le D(m+2)^{-1/2}, \qquad m\ge 0,
\ee
for some $\pi>0$ and some constant $D>0$. Then, assuming the remaining conditions of Theorem \ref{theo1} Part b), we have the estimate
\be\label{NewRandom}
E(\|u-u^{(m)}\|_a^2)^{1/2} \le 2 \|u-h\|_a + C (m+1)^{-1/2},\qquad m\ge 0,
\ee
for any $h\in {\mA^{\pi}_\infty}$, with some constant $C$ depending on $\Lambda$, $\|u\|_a$, $\|h\|_{\mA^{\pi}_\infty}$ and the constant $D$ in (\ref{PiError}).
\end{pro}

{\bf Proof}. We repeat the same steps that lead to (\ref{Erandom1}), with the following changes: The conditional expectation
$E(\|e^{(m+1)}\|_a^2\,|\,u^{(m)})$ is now computed with respect to $\pi^{(m)}$. For estimating the difference $a(e^{(m)},u)-\xi \sum_i \pi_i^{(m)} a(e^{(m)},w_i^{(m)})$, we use the $h^{(m)}\in {\mA^{\pi^{(m)}}_\infty}$ whose existence is guaranteed by Lemma \ref{lem2}, and set 
$\xi=\|h^{(m)}\|_{\mA^{\pi^{(m)}}_\infty}$. This gives
\bea
&& a(e^{(m)},u)-\xi \sum_i \pi_i^{(m)} a(e^{(m)},w_i^{(m)}) \le a(e^{(m)},u-h^{(m)}) \le \|e^{(m)}\|_a(\|u-h\|_a+ \|h-h^{(m)}\|_a)\\
%&&\qquad\qquad \le \|e^{(m)}\|_a(\|u-h\|_a+ \|h-h^{(m)}\|_a)\\
&&\qquad\qquad\qquad \le \|e^{(m)}\|_a(\|u-h\|_a + D(1+3\gL)\|h\|_{\mA^{\pi}_\infty} (m+2)^{-1/2}),
\eea
where we have substituted (\ref{Hm}) and (\ref{PiError}). Using as before the notation $\teps_m= E(\|e^{(m)}\|_a^2)^{1/2}$ and the obvious
inequality $E(\|e^{(m)}\|_a)\le \teps_m$, we arrive at the following replacement for (\ref{Erandom1}):
\be\label{Erandom3}
\teps_{m+1} \le \alpha_m \teps_m+ \bar{\alpha}_m (2\|u-h\|_a  + C_0(m+2)^{-1/2}) + \frac{\bar{\alpha}_m^2 C_1}{\alpha_m\teps_m},
\ee
where as before $\teps_m:=E(\|e^{(m)}\|_a^2)^{1/2}$. The constants are $C_0=2D(1+3\gL)\|h\|_{\mA^{\pi}_\infty}$, and
$$
2(\|u\|_a^2 + \gL^{2}\|h^{(m)}\|_{\mA^{\pi^{(m)}}_\infty}^2) \le C_1:=2(\|u\|_a^2 + 9\gL^{2}\|h\|_{\mA^{\pi}_\infty}^2),
$$
see (\ref{Hm}). This gives a recursion for 
$$
b_m:=(m+1)^{1/2}C_1^{-1/2}(\teps_m-2\|u-h\|_a),\qquad m\ge 0,
$$
similar to (\ref{E1}) but with a new term induced by the additional term $2C_0(m+1)^{-3/2}$ in (\ref{Erandom3}):
\be\label{Erandom4}
b_{m+1} \le \alpha_m^{1/2}(b_m+ \frac{1}{(m+1)b_m})+\frac{\tilde{B}}{m+2},      \qquad \tilde{B}=\frac{C_0}{C_1^{1/2}}.
\ee
This relation is again complemented by the inequality (\ref{E2}), this time with the constant $B=\|u\|_a C_1^{-1/2}\le 1/\sqrt{2}$. Since
repeating the proof of Lemma \ref{lem3} with the additional term in the right-hand side of (\ref{Erandom4}) does not represent any difficulty,
we leave it to the reader to show that $b_m\le A$, $m\ge 0$, holds for some new constant $A$ depending on $B$ and 
$\tilde{B}$. This shows (\ref{NewRandom}) with $C=AC_1^{1/2}$, and finishes our sketch of the proof of Proposition \ref{pro3}. \hfill $\Box$

\smallskip
{\bf Remark 5.} In the generality considered here, the obtained convergence rates for the error $\|u-u^{(m)}\|_a$ are not very impressive but unfortunately cannot be improved much. For greedy orderings, this issue has been addressed in \cite{DeTe,Te,Te1}. We add some comments for random orderings. Consider the very special situation of a one-dimensional subspace splitting induced by a complete orthonormal system $\mathcal{D}$ in $V$ with $a(\cdot,\cdot)=(\cdot,\cdot)$, and the problem of incremental approximation of a given $u\in V$ by linear combinations of elements from $\mathcal{D}$. I.e., if 
$$
u=\sum_{i\in\NN} c_i \psi_i 
$$
is the unique orthogonal decomposition of $u$ with respect to $\mathcal{D}$, then we have $R_iT_iu=c_i\psi_i$. Fix the discrete probability distribution $\pi>0$,
and consider the associated randomized Schwarz iterative method with updates of the form (\ref{updateLin}). It is easy to find that, due to the
orthogonality of the splitting, the best expected convergence rate is achieved for $\alpha_m=1$. In that case, we have
$$
u^{(m)} = \sum_{i\in I^{(m)}} c_i\psi_i,\qquad \|u-u^{(m)}\|^2 = \sum_{i\not\in I^{(m)}} |c_i|^2,
$$
with probability $\prod_{k=0}^{m-1} p_{i_k}$, where $\{i_k\}_{k\ge 0}$ is the random index sequence, and $I^{(m)}$ is the set of the first $m$ such indices
($I^{(m)}$ may have cardinality $\le m$, repetitions are possible). We leave it to the reader  to verify the identity
$$
E(\|u-u^{(m)}\|^2) = \sum_{i\in\NN} |c_i|^2 (1-\pi_i)^m,\qquad m\ge 1. 
$$
In the particular case considered, the formula confirms the statement of Theorem \ref{theo2} b): Since $\pi_i>0$ for all $i\in \NN$, we have $E(\|u-u^{(m)}\|^2)\to 0$
as $m\to \infty$, i.e., the expected error converges to $0$ for any $u\in V$ and any probability distribution $\pi>0$.

On the other hand, when inspecting the statement of Theorem \ref{theo1} b) for our case, we see that $u\in {\mA^{\pi}_\infty}$ is equivalent to
the inequality 
$$
|c_i|\le C\pi_i,\quad i\in \NN, \qquad C:=\|u\|_{\mA^{\pi}_\infty}<\infty.
$$
Thus, for $u\in \mA^{\pi}_\infty$ we have
\be\label{PCONS}
E(\|u-u^{(m)}\|^2) \le \|u\|_{\mA^{\pi}_\infty}^2 \sum_{i\in\NN} \pi_i^2 (1-\pi_i)^m,\qquad m\ge 1,
\ee
which is sharp in the sense that equality holds (simultaneously for all $m\ge 1$) if we set $c_i=C\pi_i$. Since
$\phi(t)=t^2(1-t)^m\le c/(m+1)^2$ for $t\in [0,1]$ for some absolute constant $c$ and all $m\ge 1$, we get
$$
\sum_{i\in\NN} \pi_i^2 (1-\pi_i)^m \le \sum_{\pi_i\le (m+1)^{-1}} \pi_i^2  + \sum_{\pi_i>(m+1)^{-1}} \frac{c}{(m+1)^2} < \frac{1+c}{m+1}.
$$
Indeed, the first sum can be estimated according to
$$
\sum_{\pi_i\le (m+1)^{-1}} \pi_i^2 \le \frac{1}{m+1}\sum_{\pi_i\le (m+1)^{-1}} \pi_i \le \frac{1}{m+1},
$$
and the second is a finite sum with $\le m$ terms. This result is in line with the bound of Theorem \ref{theo1}. 

No substantial improvement of the decay rate $O((m+1)^{-1})$  can be expected for general probability distributions: For each fixed $m$, taking $\pi$ sufficiently close to the uniform distribution
on $\{1,\ldots,m+1\}$ provides a lower bound of $c'/(m+1)$ for the right-hand side in (\ref{PCONS}) while 
choosing $\pi >0$ according to 
$$
\pi_i=c\frac1{i\log(i+1)^2},\quad i\in\NN,\qquad c:=\left(\sum_{i\in\NN} \frac1{i\log(i+1)^2}\right)^{-1},
$$
shows that, for some $u\in {\mA^{\pi}_\infty}$, lower bounds of the form $E(\|u-u^{(m)}\|^2)\ge c_\alpha (m+1)^{-\alpha}$ may hold for all $m\ge 0$ simultaneously, with any $\alpha>1$.
However, for sequences $\pi$ of the form 
$$
\pi_i = c_s i^{-(1+s)},\quad i\in\NN, \qquad c_s:=\left(\sum_{i\in\NN} i^{-(1+s)}\right)^{-1},
$$
with $s > 0$, slight improvements are possible. Note that for specific complete orthonormal systems (such as the trigonometric system in
$V=L^2(\TT)$) the classes ${\mA^{\pi}_\infty}$ for such $\pi$ have natural interpretations as $L^2$-Besov-Lipschitz spaces (in the case of our example
this would be $B_{2,\infty}^s$, and $s$ corresponds to a smoothness parameter). 
Better rates can also be concluded if we assume that $u$ belongs to a smaller class of this type with parameter $s'>s$.

Although the validity of these observations heavily relies on the assumed orthogonality of the splitting, we believe that 
especially the randomized versions should be investigated further. In particular, improved convergence rates 
for special classes of space splittings (e.g., induced by multilevel frames and sparse grid spaces) are desirable, and the potential of randomization techniques for the development of new adaptive algorithms needs to be further evaluated.

%\section{Concluding remarks}\label{sec4}

 \section*{Acknowledgement}
M. Griebel was partially supported by the project {\em EXAHD} of the 
DFG priority program 1648 {\em Software for Exascale Computing" (SPPEXA)}
and by the Sonderforschungsbereich 1060 {\em The Mathematics of Emergent Effects}
funded by the Deutsche Forschungsgemeinschaft.

\bibliographystyle{amsplain}

\end{document}